\documentclass[11pt]{article}

\usepackage{fullpage}
\usepackage{graphicx}
\usepackage{amsmath}
\usepackage{amssymb}
\usepackage{amsfonts}
\usepackage{amsthm}
\usepackage{lmodern}

\DeclareMathOperator*{\argmax}{arg\,max}

\def\ee{\mathbb{E}} % ESPERANCE
\def\bX{\mbox{\boldmath $X$}}

\def\bSigma{\mbox{\boldmath $\Sigma$}}

\begin{document} 

\begin{center}

{\bf \large Stochastic simulators based optimization by Gaussian process metamodels - Application to maintenance investments planning issues} 
\vspace{0.5cm}

{\bf Thomas BROWNE, Bertrand IOOSS, Lo\"ic LE GRATIET, J\'erome LONCHAMPT}

EDF R\&D, 6 Quai Watier, F-78401 Chatou, France 
\end{center}
 
\vspace{0.5cm}

\abstract{
This paper deals with the construction of a metamodel (i.e. a simplified mathematical model) for a stochastic computer code (also called stochastic numerical model or stochastic simulator), where stochastic means that the code maps the realization of a random variable. 
The goal is to get, for a given model input, the main information about the output probability distribution by using this metamodel and without running the computer code. 
In practical applications, such a metamodel enables one to have estimations of every possible random variable properties, such as the expectation, the probability of exceeding a threshold or any quantile. 
The present work is concentrated on the emulation of the quantile function of the stochastic simulator by interpolating well chosen basis function and metamodeling their coefficients (using the Gaussian process metamodel).
This quantile function metamodel is then used to treat a simple optimization strategy maintenance problem using a stochastic code, in order to optimize the quantile of an economic indicator. 
Using the Gaussian process framework, an adaptive design method (called QFEI) is defined by extending in our case the well known EGO algorithm. 
This allows to obtain an ``optimal'' solution using a small number of simulator runs.
}
 
\vspace{0.5cm}
\noindent {\bf Keywords:} Adaptive Design, Asset management, Computer experiments, Expected Improvement, Gaussian process, Optimization, Stochastic simulator, Uncertainty.

%\clearpage
%%%%%%%%%%%%%%%%%%%%%%%%%%%%%%%%%%%%%%%%%%%%%%%%%%%%%%%%%%%%% 

\section{Introduction}

EDF looks for assessing and optimizing its strategic investments decisions for its electricity production assets by using probabilistic and optimization methods of ``cost of maintenance strategies''. 
In order to quantify the technical and economic impact of a candidate maintenance strategy, some economic indicators are evaluated by Monte Carlo simulations using the VME software developed by EDF R\&D (Lonchamp and Fessart \cite{lonchamptfressart:2013}). 
The major output result of the Monte Carlo simulation process from VME is the probability distribution of the Net Present Value (NPV) associated to the maintenance strategy. 
From this distribution, some indicators, such as the NPV mean, the NPV standard deviation or the regret investment probability ($\mathbb{P}(NPV < 0)$) can easily be derived.

Once these indicators have been obtained, one is interested in optimizing the strategy, for instance by determining the optimal investments dates leading to the highest mean NPV and the lowest regret investment probability. 
Due to the discrete nature of the events to be optimized, the optimization method is actually based on genetic algorithms.
However, during the optimization process, one of the main issues is the computational cost of the stochastic simulator to optimize, which leads to methods requiring a minimal number of simulator runs (Dellino and Meloni \cite{delmel15}).
Genetic algorithms require often several thousands of simulator runs and, in some cases, are no more applicable.

The solution investigated in this study is a first attempt to break the computational cost of this problem by the way of using a metamodel instead of the simulator within the mathematical optimization algorithm. 
From a first set of simulator runs (called the learning sample and coming from a specific design of experiments), a metamodel consists in approximating the simulator outputs by a mathematical model (Fang et al. \cite{fanli06}). This metamodel can then be used to predict the simulator outputs for other input configurations. 

Many metamodeling techniques are available in the computer experiments literature (Fang et al. \cite{fanli06}).
Formally, the function $G$ representing the computer model is defined as 
 \begin{equation}
\begin{array}{rccl}
 G :& E & \rightarrow & \mathbb{R} \\
&x &\mapsto & G(x) 
\end{array}
\end{equation}
where $E\subset \mathbb{R}^d$ ($d\in \mathbb{N}^*$) is the space of the input variables of the computer model.
Metamodeling consists in the construction of a statistical estimator $\widehat{G}$ from an initial sample of $G$ values corresponding to a learning set $\chi$ with $\chi \subset E$ and $\# \chi < \infty$ (with $\#$ the cardinal operator).

However, the conventional methods are not suitable in the present framework because of the stochastic nature of the simulator: the output of interest is not a single scalar variable but a full probability density function (or a cumulative distribution function, or a quantile function).
The computer code $G$ is therefore stochastic:
\begin{equation}
\begin{array}{rccl}
G \, : \,& E\times\Omega & \rightarrow & \mathbb{R} \\
& (x,\omega) & \mapsto & G(x,\omega) 
\end{array}
\end{equation}
where $\Omega$ denotes the probability space. 
In the framework of robust design, robust optimization and sensitivity analysis, previous works with stochastic simulator consist mainly in approximating the mean and variance of the stochastic output (Bursztyn and Steinberg \cite{burste06}, Kleijnen \cite{kle08}, Ankenman et al \cite{anknel10}, Marrel et al. \cite{marioo12}, Dellino and Meloni \cite{delmel15}).
Other specific characteristics of the distribution of the random variable $G(x)$ (as a quantile) can also be modeled to obtain, in some cases, more efficient algorithms (Picheny et al. \cite{picgin13}).

In this paper, as first proposed by Reich et al. \cite{reikal12} (who used a simple Gaussian mixture model), we are interested in a metamodel which predicts the full distribution of the random variable $G(x)$, $\forall x \in E$. 
We then define the following deterministic function $f$:
\begin{equation}
\begin{array}{rccl}
f \, : \,& E & \rightarrow & \mathcal{F} \\
& x & \mapsto & f_x \quad \text{density of the r.v. $G(x)$}  \\
\end{array}
\end{equation}
with $\mathcal{F}$ the family of probability density functions:
 \begin{equation}  \mathcal{F} = \left\{ g \in L^1(\mathbb{R}), \; \text{positive, measurable}, \int_{ \mathbb{R} } g = 1 \right\} .  \end{equation} 
For a point $x \in E$, building $f_x$ with the kernel method requires a large number $N_{\mbox{\scriptsize MC}}$ realizations of $G(x)$.
A recent work (Moutoussamy et al. \cite{mounan14}) has proposed kernel-regression based algorithms to build an estimator $\widehat{f}$ of the output densities, under the constraint that each call to $f$ is cpu-time expensive.
Their result stands as the starting point for the work presented in this paper (based on a Master internship report, see Browne \cite{bro14}).

The next section briefly presents the VME application case.
In the third section, we propose to use the Gaussian process metamodel and develop an algorithm to emulate the quantile function instead of the probability density function. 
In the fourth section, this metamodel is used to treat the VME case, in order to optimize a quantile. 
Using the Gaussian process metamodel framework, an adaptive design method is also proposed, allowing to obtain an ``optimal'' solution using a small number of VME simulator runs.
A conclusion synthesizes the main results of this paper.

%*********************************************************
\section{The VME application case}

On an industrial facilities, we are interested by the replacement cost of four components in function of the date of their maintenance (in year) (see Lonchamp and Fessart \cite{lonchamptfressart:2012} for more details).
To these four inputs, one additional input is the date (in year) of recovering a new component.
This input space is denoted $F$  : 
	\begin{equation}  F=\left(\bigotimes _{i=1}^4 \left\lbrace 41,42,...,50 \right\rbrace  \right) \times  \left\lbrace 11,12,...,20 \right\rbrace .  \end{equation} 
$F$ is therefore a discrete set ($\# F = 10^5 $ ).
We have
	\begin{equation} \label{eq:densmodel}
\begin{array}{rccl}
G \, : \,& F & \rightarrow & \mathbb{R} \\
& x=(x_1,x_2,x_3,x_4,x_5) & \mapsto & \text{NPV}(x) , \\
\\
f \, : \,& F & \rightarrow & \mathcal{F} \\
& x=(x_1,x_2,x_3,x_4,x_5) & \mapsto & f_x \quad \text{(density of NPV}(x))  .\\
\end{array}
 \end{equation} 

\begin{figure}[!ht]
\begin{center}
\includegraphics[width=8cm,height=6cm]{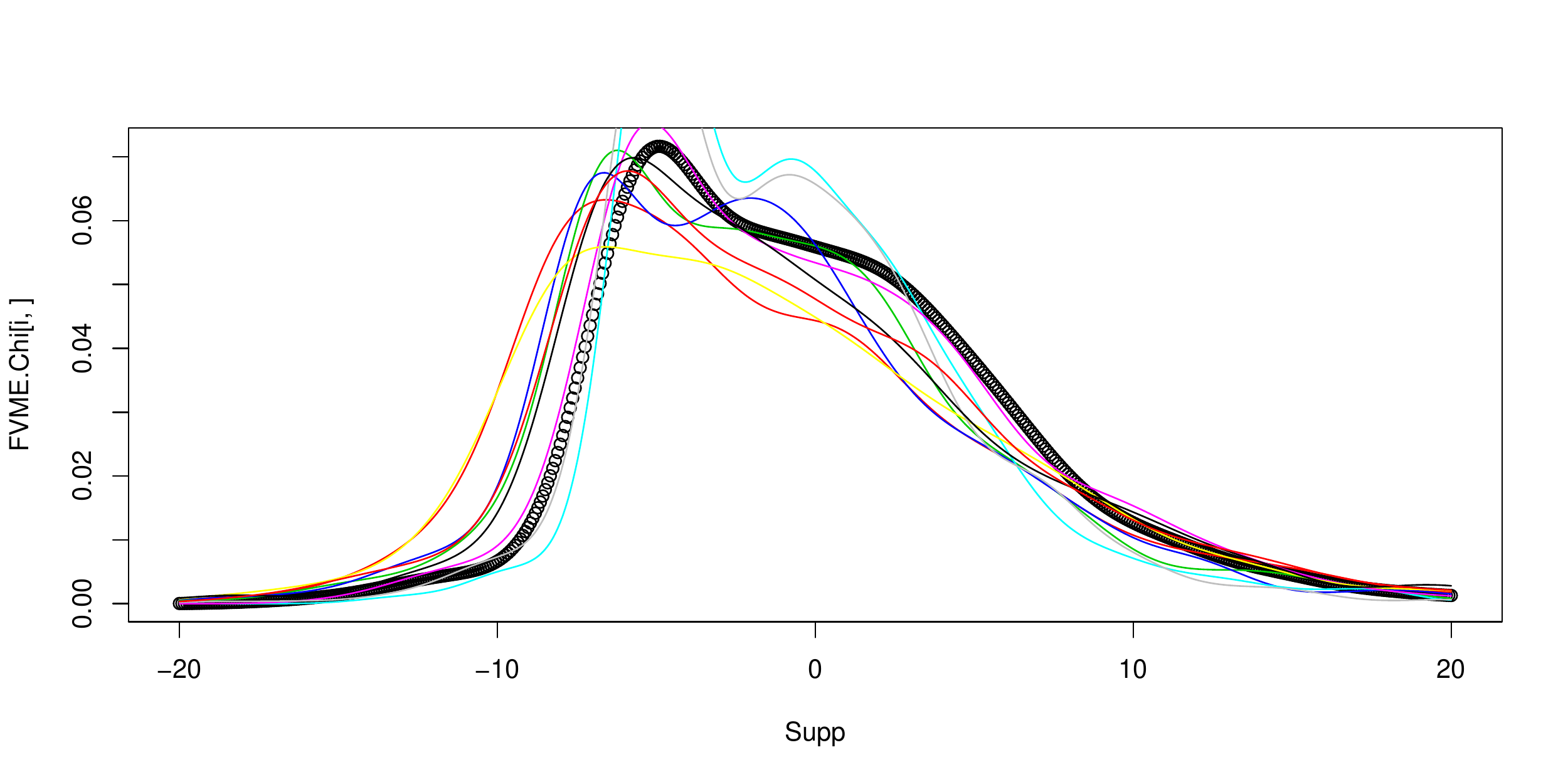}
\caption{$10$ ouput probability densities of the $VME$ code (randomly sampled).} 
\label{Opt_dens}
\end{center}
\end{figure}

Figure \ref{Opt_dens} provides examples of the output density of VME.
The $10$ input values inside $F$ have been randomly chosen.
It shows that there is a small variability between the curves.

The optimization process consists in finding the point of $F$ which gives the ``maximal NPV'' value. 
As NPV$(x)$ is a random variable, we have to summarize its distribution by a deterministic operator $H$, for example:
	\begin{equation}   
 H(g)=\ee(g) \; \forall g \in \mathcal{F}
 \end{equation}   
or
	\begin{equation}   
 H(g)=q_g(\alpha)\; \forall g \in \mathcal{F}
 \end{equation}   
with $q_g(\alpha)$ the $\alpha$-quantile of $g$.
Our VME-optimization problem turns then to the determination of
\begin{equation}  x^* := \argmax_{x \in F} H(f_x)  . \end{equation} 

However, several difficulties occur:
\begin{itemize}
\item VME is a cpu-time costly simulator and the size of the set $F$ is large. 
Computing $(f_x)_{x \in F}$, needing $N_{\mbox{\scriptsize MC}} \times \# F$ simulator calls (where $N_{\mbox{\scriptsize MC}}$ is worth several thousands), is therefore impossible.
Our solution is to restrict the VME calls to a learning set $\chi \subset F$ (with $\# \chi = 200$ in our application case), randomly chosen inside $F$.
We will then have $(f_x)_{x \in \chi}$.

\item
Our metamodel, that we will denote $\left(\hat{\hat{f_x}}\right)_{x  \in F} $, cannot be computed on $F$ due to its large size.
A new set $E$ is therefore defined, with $ E \subset F$ and $\#E = 2000$.
We will then have $\left(\hat{\hat{f_x}} \right)_{x  \in E}$.
In this work, we limit our study to this restricted space $E$ instead of the full space $F$.
Other work will be performed to extend our methodology to the $F$ study.
\end{itemize}

%*********************************************************
\section{Gaussian process metamodel of a stochastic simulator}

%%%%%%%%%%%%%%%%%%%%%%%%%%
\subsection{Basics on the Gaussian process model}

Let us consider $n$ evaluations of a deterministic computer code. 
Each evaluation $Y(x) \in \mathbb{R}$ of a simulator scalar output comes from a $d$-dimensional input vector $x = (x_1,\ldots,x_d) \in E$, where $E$ is a bounded domain of $\mathbb{R}^d$. 
The $n$ points corresponding to the code runs are called the experimental design and are denoted as $\bX_s = ( x^{(1)},\ldots,x^{(n)})$. 
This learning sample comes from the learning set, previously noted $\chi$ ($n=\#\chi$).
The outputs will be denoted as $Y_s = (y^{(1)},\ldots,y^{(n)})$ with $y^{(i)} = G(x^{(i)}) \; \forall \; i=1..n$.
Gaussian process modeling (Sacks et al. \cite{sacwel89}), also called kriging model (Stein \cite{ste99}), treats the simulator deterministic response $G(x)$ as a realization of a random function $Y(x)$, including a regression part and a centered stochastic process:
\begin{equation}\label{Gpmodel}
 Y ( x ) = h ( x ) + Z ( x) .
\end{equation}

The deterministic function $h(x)$ provides the mean approximation of the computer code. We can use for example a one-degree polynomial model:
\begin{equation}  h ( x ) = \beta_0 + \sum_{j = 1}^d \beta_j x_j \;,  \end{equation} 
where $\beta = [ \beta_0, \ldots, \beta_d ]^t $ is the regression parameter vector.
The stochastic part $Z(x)$ is a Gaussian centered stationary process fully characterized by its covariance function:
$\mbox{Cov} ( Z ( x ), Z ( u ) ) = \sigma^2 K_{\theta} ( x - u ),$
where $\sigma^2$ denotes the variance of $Z$, $K_{\theta}$ is the correlation function and $\theta \in \mathbb{R}^d$ is the vector of correlation  hyperparameters.
This structure allows to provide interpolation and spatial correlation properties.
Several parametric families of correlation functions can be chosen (Stein \cite{ste99}).

If a new point $x^{\ast} = (x^{\ast}_1,\ldots,x^{\ast}_d) \in E$ is considered, we obtain the predictor and variance formulas for the scalar output $Y(x^*)$:
 \begin{eqnarray}
 \displaystyle \ee [  Y(x^{*})|Y_s] = h(x^{*}) +  k(x^{*})  ^t \bSigma_s^{-1} (Y_s - h(\bX_s)) \;, \label{eq_esperance} \\ 
 \displaystyle MSE(x^{*}) = \mbox{Var}[ Y(x^{*})|Y_s] = \sigma^2  -  k(x^{*}) ^t  \bSigma_s^{-1} k(x^{*}) \;, \label{eq_variance}
\end{eqnarray}
with
  \begin{equation}  
\begin{array}{lll}
k(x^{*}) & = &[\mbox{Cov}(y^{(1)},Y(x^{*})), \ldots, \mbox{Cov}(y^{(n)},Y(x^{*}))  ] ^t  \\
& = & \sigma^2  [  K_{\theta} (x^{(1)},x^{*}), \ldots, K_{\theta} (x^{(n)},x^{*}) )  ]^t  
 \end{array} 
 \end{equation} 
and the covariance matrix
\begin{equation}  \bSigma_s = \sigma^2 \left( K_{\theta} \left( x^{(i)} - x^{(j)} \right)  \right)_{i=1\ldots n, j = 1\ldots n} \;. \end{equation} 
The conditional mean (Eq. (\ref{eq_esperance})) is  used as a predictor. The variance formula (Eq. (\ref{eq_variance})) corresponds to the mean squared error (MSE) of this predictor and is also known as the kriging variance. This analytical formula for MSE gives a local indicator of the prediction accuracy. 
More generally, Gaussian process model defines a Gaussian distribution for the output variable at any arbitrary new point. 
This distribution formula can be used for uncertainty and sensitivity analysis (Marrel et al. \cite{marioo08}, Le Gratiet et al. \cite{legcan14}). 
Regression and correlation parameters $\beta$, $\sigma$ and $\theta$ are ordinarily estimated by maximizing likelihood functions (Fang et al. \cite{fanli06}).

%%%%%%%%%%%%%%%%%%%%%%%%%%
\subsection{Emulation of the simulator quantile function}

%%%%%
\subsubsection{General principles}

\begin{figure}[!ht]
\begin{center}
\includegraphics[width=8cm,height=6cm]{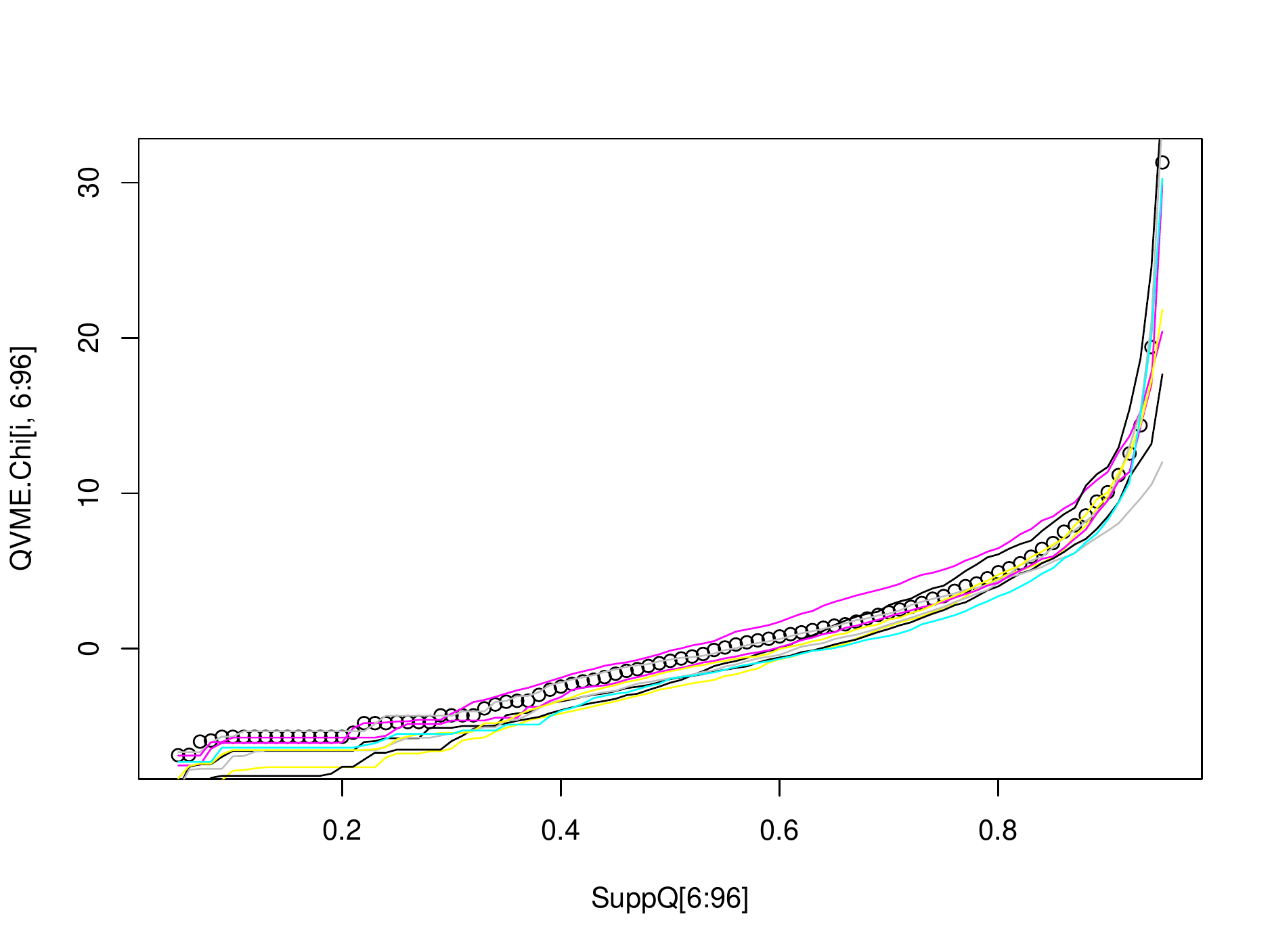}
\caption{Quantile functions $Q$ for 10 points of $E$ (randomly chosen).} \label{plot5}
\end{center}
\end{figure}

In our VME-optimization problem, we are especially interested by several quantiles (for example at the order $1\%$, $5\%$, $50\%$, $95\%$, $99\%$) rather than statistical moments.
In Moutoussamy et al. \cite{mounan14} and Browne \cite{bro14}, quantile prediction with density-based emulator has shown some deficiencies.
Therefore, instead of studying Eq. (\ref{eq:densmodel}), we turn our modeling problem to
\begin{equation} 
\begin{array}{rccl}
G \, : \,& E & \rightarrow & \mathbb{R} \\
& x=(x_1,x_2,x_3,x_4,x_5) & \mapsto & \text{NPV}(x) , \\
\\
Q \, : \,& E & \rightarrow & \mathcal{Q} \\
& x & \mapsto & Q_x \quad \text{quantile function of NPV$(x)$}  
\end{array}
 \end{equation} 
where $\mathcal{Q}$ is the space of increasing functions defined on $]0,1[$, with values in $[a,b]$ (which is the support of the NPV output). 
For $x \in E$, a quantile function is defined by :
\begin{equation}  
\forall p \in ]0,1[, \;\; Q_x(p) = t \in [a,b] \quad \mbox{such as} \quad \int_a^t f_x(\varepsilon)d\varepsilon=p   .
\end{equation}  
For the same points than in Figure \ref{Opt_dens}, Figure \ref{plot5} shows the $10$ quantile function ouputs $Q$ which present a rather low variability.

We consider the learning set $\chi$ ($n=\#\chi$) and $N_{\mbox{\scriptsize MC}} \times n$ $G$-simulator calls in order to obtain $\left( \tilde{Q}_x^{N_{\mbox{\tiny MC}}} \right)_{x \in \chi}$, the empirical quantile functions of $\left( \mbox{NPV}(x) \right)_{x \in \chi}$. 
In this work, we will use $N_{\mbox{\scriptsize MC}} = 10^4$, which is sufficiently large to obtain a precise estimator of $Q_x$ with $\tilde{Q}_x^{N_{\mbox{\tiny MC}}}$.
Therefore, we neglect this Monte Carlo error.
In the following, we simplify the notations by replacing $\tilde{Q}_x^{N_{\mbox{\tiny MC}}}$ by $Q_x$.

The approach we adopt is similar to the one used in metamodeling a functional output of a deterministic simulator (Bayarri et al., \cite{bayber07}, Marrel et. al. \cite{marioo11}).
The first step consists in finding a low-dimensional functional basis in order to reduce the output dimension by projection, while the second step consists in emulating the coefficients of the basis functions.
However, in our case, due to the nature of the functional outputs (quantile functions), some particularities will arise. 

%%%%%
\subsubsection{Projection of $Q_x$ by the Modified Magic Points (MMP) algorithm}

Adapted from the Magic Points algorithm (Maday et al. \cite{madngu07}) for probability density functions, the MMP algorithm has been proposed in Moutoussamy et al. \cite{mounan14}.
It is a greedy algorithm that builds interpolator (as a linear combination of basis functions) for a set of functions by iteratively picking a basis function in the learning sample output set and a set of interpolation points.
At each step $j \in \{2, \ldots , q\}$ of the construction of the functional basis, one picks the element of the learning sample output set that maximizes the gap ($L^2$ distance) between this element and the interpolator which used the previous $j-1$ functions of the basis.
The total number $q$ of functions is chosen with respect to a convergence criterion.
Mathematical details will not be provided in the present paper.

In this paper, we apply the MMP algorithm on quantile functions.
The first step consists in a projection of $\left( Q_x \right)_{x \in \chi}$:
\begin{equation} \hat{Q}_x = \sum_{j=1}^q \psi_j(x) R_j \; \forall x \in \chi 
 \end{equation} 
where $\psi=(\psi_1,\ldots,\psi_q)$ (the coefficients) and $R=\left( R_1,...,R_q \right)$ (the quantile functions of the basis) are determined by MMP.
We need to restrict the solutions to the following constrained space:
\begin{equation}  
C = \left\lbrace \psi \in \mathbb{R}^q, \quad \psi_1,...,\psi_q \geq 0 \right\rbrace  
\end{equation} 
in order to ensure the monotonic increase of $\hat{Q}_x$.

\begin{figure}[!ht]
\begin{center}
\includegraphics[width=8cm,height=6cm]{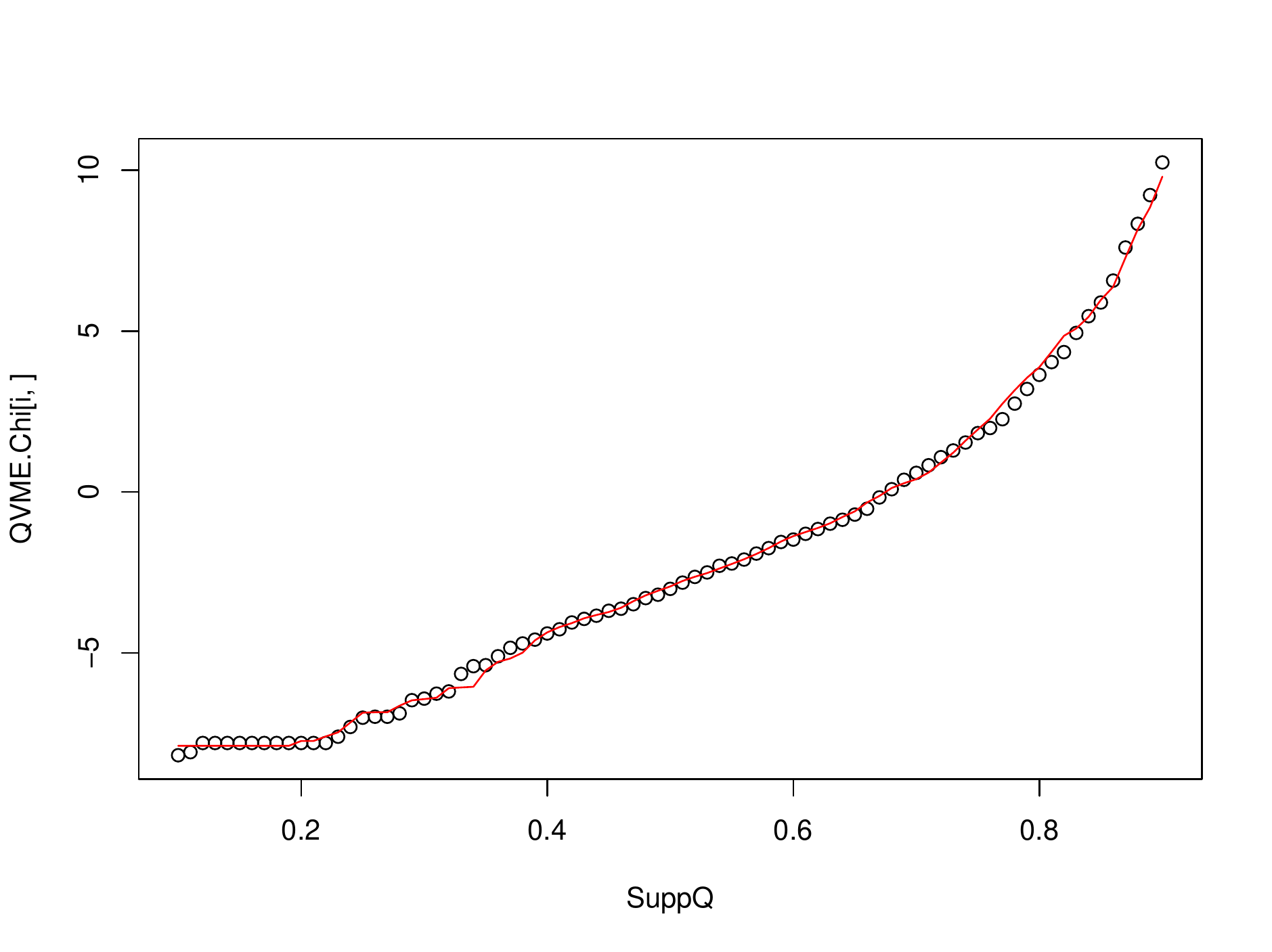}
\caption{For one point $x \in \chi$, $\hat{Q}_x$ (red line) and $Q_x$ (black points).} \label{plot6}
\end{center}
\end{figure}

In our VME application, the choice $q=5$ has shown sufficient approximation capabilities.
For one example of quantile function output, a small relative $L^2$-error ($0.2\%$) between the observed quantile function and the projected quantile function is obtained.
Figure \ref{plot6} confirms also the relevance of the MMP method.

%%%%%
\subsubsection{Gaussian process metamodeling of the basis coefficients}\label{sec:metamodelcoef}
	
Estimations of the $\psi(x) = (\psi_1(x), \ldots, \psi_q(x))$ ($x \in E$) coefficient vector will be performed with $q$ Gaussian process metamodels in order to build $\hat{\hat{Q}}_x$:
\begin{equation}\label{eq:metamodel}
  \hat{\hat{Q}}_x = \sum_{j=1}^q \hat{\psi}_j(x)R_j  \; \forall x \in E . 
	\end{equation} 
To ensure that $\hat{\psi}_j \in C$ $(j=1\ldots q$), we use a logarithmic transformation:
\begin{equation} 
\begin{array}{rccl}
{\cal T}_1 \, : \,& C & \rightarrow & \mathbb{R}^q \\
& \psi & \mapsto & \left(\log(\psi_1+1),...,\log(\psi_q+1) \right) \\
\end{array}
 \end{equation} 
 and its inverse transformation:
 \begin{equation} 
\begin{array}{rccl}
{\cal T}_2 \, : \,& \mathbb{R}^q & \rightarrow & C \\
& \varphi & \mapsto & \left(\exp(\phi_1)-1,...,\exp(\phi_q)-1 \right). \\
\end{array}
 \end{equation} 
We then compute $\phi(x) := {\cal T}_1(\psi(x))  \; \forall x \in \chi $ and suppose that $\phi$ is a Gaussian process realization with $q$ independent margins. 
$\phi$ is estimated by 
\begin{equation}
\hat{\phi}(x):= \ee[Y_x \mid Y_s ] \; \forall x \in E 
\end{equation}
 with $Y_s$ the learning sample output.
We obtain 
\begin{equation}
\hat{\psi}(x):= {\cal T}_2 (\hat{\phi}(x)) \; \forall x \in E 
\end{equation}
and Eq. (\ref{eq:metamodel}) can be applied as our metamodel predictor of the quantile function.

\begin{figure}[!ht]
\begin{center}
\includegraphics[width=8cm,height=6cm]{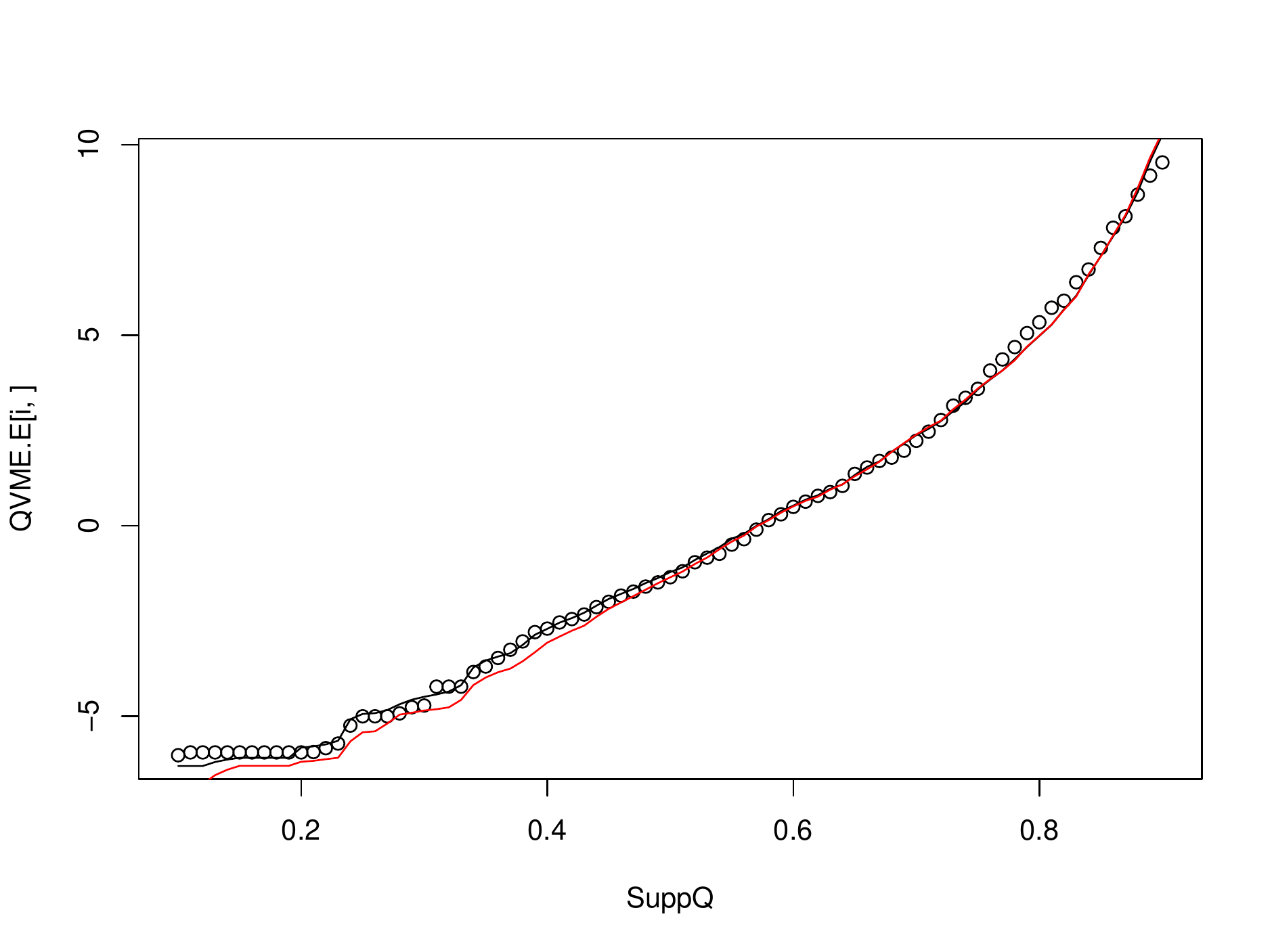}
\caption{For one point $x \in \chi$, $\hat{\hat{Q}}_x$ (red line) and $Q_x$ (black points).}\label{plot7}
\end{center}
\end{figure}

In our VME application, we have built the metamodel on the set $E$ (with the choice $q=5$).
For one example of quantile function output, a small relative $L^2$-error ($2.8\%$) between the observed quantile function and the emulated quantile function is obtained.
Figure \ref{plot7} confirms also the relevance of the metamodeling method.

%*********************************************************
\section{Application to an optimization problem}

%%%%%%%%%%%%%%%%%%%%%%%%%%
\subsection{Direct optimization on the metamodel}

We now apply our quantile function metamodel with a quantile-based objective function
\begin{equation}   
\begin{array}{rccl}
H  \, : \,& \mathcal{Q} & \rightarrow & \mathbb{R} \\
& q & \mapsto & q(p) \\
\end{array}  
 \end{equation}   
with $p\in ]0,1[$.
We look for
\begin{equation}  x^* := \argmax_{x \in F} Q_x(p) \end{equation} 
but have only access to
\begin{equation}  \hat{\hat{x}}^* := \argmax_{x \in F} \hat{\hat{Q}}_x(p) . \end{equation} 
%On aimerait avoir $\hat{\hat{x}}^*\simeq x^* $ mais pour cela il faut absolument que l'on ait :
%\begin{equation}  
 %\quad H(\hat{\hat{Q}}_x) \simeq H(Q_x) \; \forall x \in F.
  %\end{equation} 
	
We study the relative error of $H(\hat{\hat{Q}})$ on $E$ by computing
%Si elle est déjà trop grande alors on ne pourra chercher à étendre le métamodèle sur tout $F$ et notre méthode de résolution sera inefficace.
\begin{equation}  
err = \frac{1}{\max_{x \in E} \left( Q_x(p) \right) -\min_{x \in E} \left( Q_x(p) \right) } \times \left( \sum_{x \in E} \mid Q_x(p) - \hat{\hat{Q}}_x(p) \mid \right) .
\end{equation} 
As an example, for $p=0.5$ (median estimation), we find
\begin{equation}   
\begin{array}{lcr}
 \max_{x \in E} (  Q_x(p) ) =0.82,  & \max_{x \in E} ( \hat{\hat{Q}}_x(p) )=0.42, & err = 5.4\% . 
\end{array}  
 \end{equation}   
If we define $y= \argmax_{x \in E} \hat{\hat{Q}}_x(p)$ the best point from the metamodel, we obtain $Q_y(p)=0.29$ while $\max_{x\in \chi}Q_x(p)=0.35$. 
The exploration of $E$ by our metamodel does not bring any information. 
We have observed the same result by repeating $100$ times the experiments (changing the initial design).
It means that the punctual errors on the quantile function metamodel are too large for this optimization algorithm.
In fact, the basis functions $R_1,...,R_5$ that the MMP algorithm has chosen on $\chi$ are not able to represent the extreme parts of the quantile functions of $E$.

As a conclusion of these tests, the quantile function metamodel cannot be directly applied to solve the optimization problem.
In the next section, we propose an adaptive algorithm which consists in sequentially adding simulation points in order to capture interesting quantile functions to be added in our functional basis.

%************************
\subsection{QFEI: An adaptive optimization algorithm }

After the choice of $\chi$, $E$ and the families $\left( Q_x \right)_{x \in \chi}$, $( \hat{Q}_x)_{x \in \chi}$ and $( \hat{\hat{Q}}_x)_{x \in E}$, our new algorithm will propose to perform new interesting (for our specific problem) calls to the VME simulator on $E$ (outside of $\chi$).
With the Gaussian process metamodel, which provides a predictor and its uncertainty bounds, this is a classical approach used for example in black-box optimization problem (Jones et al. \cite{jonsch98}) and rare event estimation (Bect et al. \cite{becgin12}).
The goal is to provide some algorithms which mix global space exploration and local optimization.

Our algorithm is based on the so-called EGO (Efficient Global Optimization) algorithm (Jones et al. \cite{jonsch98}) which uses the Expected Improvement (EI) criterion to optimize a deterministic simulator.
Our case is different as we want to maximize
\begin{equation}   
\begin{array}{rccl}
H  \, : \,& E & \rightarrow & \mathbb{R} \\
& x & \mapsto & Q_x(p) \quad \mbox{$p$-quantile of NPV$(x)$} .
\end{array}  
 \end{equation}
We will then propose a new algorithm called the QFEI (for Quantile Function Expected Improvement) algorithm.

As previously, we use the set $E \subset F$ with $\#E=5000$ ($E$ is a random sample in $F$), the initial learning set $\chi \subset F$ with $\# \chi=200$ (initial design of experiment), $\left( Q_x \right)_{x \in \chi}$, $( \hat{Q}_x)_{x \in \chi}$ and $( \hat{\hat{Q}}_x)_{x \in E}$. 
We denote $D$ the current learning set (the initial learning set increased with additional points coming from QFEI).
As Gaussianity will be needed on the components of $\psi$, we did not performed a logarithmic transformation as in Section \ref{sec:metamodelcoef}.
In our case, it has not implied negative consequences.

We apply the Gaussian process metamodeling on the $q$ independent components $\psi_1,...,\psi_q$:
\begin{equation}
 \psi_j(x) \sim \mathcal{N}\left( \hat{\psi}_j(x), MSE_j(x) \right) \; \forall j \in \{1,...,q\} \; \forall x \in E.
  \end{equation}  
As $\hat{Q}_x(p) = \sum_{j=1}^q \psi_j(x) R_j(p) $, we have
\begin{equation}
\hat{Q}_x(p) \sim \mathcal{N} \left( \hat{\hat{Q}}_x(p), \sum_{j=1}^q R_j(p)^2 MSE_j(x)  \right)  \; \forall x \in E. 
\end{equation} 
Then $\hat{Q}_{x}(p)$ is a realization of the underlying Gaussian process $U_{x}=\sum_{j=1}^q \psi_j (x) R_j(p)$ with
\begin{equation} 
\begin{array}{lll}
U_D & := \left( U_x \right) _{x \in D} &,\\
\hat{U}_{x} & := \ee[U_{x} \mid U_D] &\forall x \in E ,\\
\sigma_{U|D}^2(x) & := \mbox{Var}[U_{x} \mid U_D] &\forall x \in E .
\end{array} \end{equation} 
The conditional mean and variance of $U_x$ are directly obtained from the $q$ Gaussian process metamodels of the $\psi$ coefficients.

At present, we propose to use the following improvement random function:
\begin{equation}   
\begin{array}{rccl}
I  \, : \,& E & \rightarrow & \mathbb{R} \\
& x & \mapsto & \left(U_x-\max\left(U_D \right)\right)^+ . 
\end{array}  
 \end{equation}   
In our adaptive design, finding a new point consists in solving:
\begin{equation}  x_{new} := \argmax_{x \in E} \ee[I(x)] .  \end{equation} 
Added points are those which have more chance to improve the current optimum. 
The expectation of the improvement function writes (the simple proof is given in Browne \cite{bro14}):
\begin{equation}  
\quad \ee[I(x)]=\sigma_{U|D}(x) \left( u(x) \phi(u(x)) + \varphi(u(x)) \right) \; \forall x \in E ,\quad \mbox{with} \quad u(x)= \frac{\hat{U}_x-\max(U_D)}{\sigma_{U|D}(x)}   
\end{equation} 
where $\varphi$ and $\phi$ correspond respectively to the density and distribution functions of the reduced centered Gaussian law. 

In practice, several iterations of this algorithm are performed, allowing to complete the experimental design $D$.
At each iteration, a new projection functional basis is computed and the $q$ Gaussian process metamodels are re-estimated.
The stopping criterion of the QFEI algorithm can be a maximal number of iterations or a stabilization criterion on the obtained solutions.
No garantee on convergence of the algorithm can be given.
In conclusion, this algorithm provides the following estimation of the optimal point $x^*$:
\begin{equation} \hat{x}^* := \argmax_{x \in D}(U_D) ,  \end{equation} 
	
In our application case, we have performed all the simulations in order to know $\left( Q_x \right) _{x \in E}$, therefore the solution $x^*$.
Our first objective is to test our proposed algorithm for $p=0.4$ which has the following solution:
\begin{equation}  \left\lbrace
\begin{array}{ll}
x^* & = \left( 41,47,48,45,18 \right) \\
Q_{x^*}(p) & = -1.72 .
\end{array}
\right. \end{equation} 
We have also computed
\begin{equation}  \frac{1}{\# E} \sum_{x \in E} Q_x(p)=-3.15, \quad \mbox{Var} \left( \left(Q_x(p)\right)_{x \in E} \right)=0.59 . \end{equation} 
We start with $D := \chi$ and we obtain
\begin{equation} \max_{x \in \chi} \left( Q_x \right) = -1.95  \end{equation} 
After $50$ iterations of the QFEI algorithm, we obtain:
\begin{equation}  \left\lbrace
\begin{array}{ll}
\hat{x}^* & = \left( 41,47,45,46,19 \right) \\
Q_{\hat{x}^*}(p) & = -1.74 .
\end{array}
\right. \end{equation} 
We observe that $\hat{x}^* \simeq x^*$ and $Q_{\hat{x}^*}(p) \simeq Q_{x^*}(p)$ which is a first confirmation of the relevance of our method. 
With respect to the initial design solution, the QFEI has allowed to obtain a strong improvement of the proposed solution.
$50$ repetitions of this experiment (changing the initial design) has also proved the robustness of QFEI.
The obtained solution is always one of the five best points on $E$.

QFEI algorithm seems promising but a lot of tests remain to perform and will be pursued in future works: changing $p$ (in particular testing extremal cases), increasing the size of $E$, increasing the dimension $d$ of the inputs, \ldots

%%%%%%%%%%%%%%%%%%%%%%%%%%%%%%%%%%%%%%%%%%%%%%%%%%%%%%%%%%%%% 
\section{Conclusion}

In this paper, we have proposed to build a metamodel of a stochastic simulator using the following key points: 
\begin{enumerate}
\item Emulation of the quantile function which proves better efficiency for our problem than the emulation of the probability density function; 
\item Decomposition of the quantile function in a sum of the quantile functions coming from the learning sample outputs; 
\item Selection of the most representative quantile functions of this decomposition using an adaptive choice algorithm (called the MMP algorithm) in order to have a small number of terms in the decomposition; 
\item Emulation of each coefficient of this decomposition by a Gaussian process metamodel, by taking into account constraints ensuring that a quantile function is built.
\end{enumerate}

The metamodel is then used to treat a simple optimization strategy maintenance problem using a stochastic simulator (VME), in order to optimize an output (NPV) quantile. 
Using the Gaussian process metamodel framework and extending the EI criterion to quantile function, the adaptive QFEI algorithm has been proposed. 
In our example, it allows to obtain an ``optimal'' solution using a small number of VME simulator runs.

This work is just a first attempt and needs to be continued in several directions:
\begin{itemize}
\item Consideration of a variable $N_{\mbox{\scriptsize MC}}$ whose decrease could help to fight against the computational cost of the stochastic simulator,
\item Improvement of the initial learning sample choice by replacing the random sample by a space filling design (Fang et al. \cite{fanli06}),
\item Algorithmic improvements to counter the cost of the metamodel evaluations and to increase the size of the study set $E$, 
\item Multi-objective optimization (several quantiles to be optimized) in order to  take advantage of our powerful quantile function emulator,
\item Application to more complex real cases,
\item Consideration of a robust optimization problem where environmental input variables of the simulator has not to be optimized but just create an additional uncertainty on the output.
\end{itemize}

%%%%%%%%%%%%%%%%%%%%%%%%%%%%%%%%%%%%%%%%%%%%%%%%%%%%%%%%%%%%% 
\section{Acknowledgments}

We are grateful to Emmanuel Remy for his comments on this paper.

% ***********************************************************

\bibliographystyle{plain}
%\bibliography{bibliorapstage,BIBLIOGRAPHIE/books,BIBLIOGRAPHIE/articles,BIBLIOGRAPHIE/inbooks,BIBLIOGRAPHIE/NT,BIBLIOGRAPHIE/thesis,BIBLIOGRAPHIE/proceedings,BIBLIOGRAPHIE/rapports,BIBLIOGRAPHIE/soft}

\end{document}